\documentclass[11 pt]{amsart}
\usepackage{amssymb}
\usepackage{amscd}
\usepackage{epsfig}
\newtheorem{theorem}{Theorem}
\newtheorem{corollary}[theorem]{Corollary}
\newtheorem{lemma}[theorem]{Lemma}

\theoremstyle{definition}

\numberwithin{equation}{section}

\newcommand{\bz}{\mathbb{Z}}
\newcommand{\br}{\mathbb{R}}

\newcommand{\bh}{\mathbb{H}}

\newcommand{\cc}{\mathcal{C}}

\newcommand{\cf}{\mathcal{F}}
\newcommand{\bx}{\mathbb{X}}

\newcommand{\hk}{\hookrightarrow}

\newcommand{\med}{\medskip}
\newcommand{\la}{\longrightarrow}
\newcommand{\bfl}{\begin{flushleft}}
\newcommand{\efl}{\end{flushleft}}

\newcommand{\ltm}{LM^{-TM}}
\newcommand{\mtm}{M^{-TM}}
\newcommand{\tbx}{(\mathbb{L}_M)}
\newcommand{\tmu}{\tilde\mu}
\newcommand{\lmm}{LM \times_M LM}

\title{A homotopy theoretic realization of String Topology}
\author[R.L. Cohen]{Ralph L. Cohen}
 \address{Dept. of Mathematics \\
Stanford University\\
Stanford, California 94305}
 \email[Cohen]
{ralph@math.stanford.edu}
\thanks{The first author was partially supported by a grant from the NSF
} 
\author[J.D.S Jones]{John D.S. Jones}
\address{Department of Mathematics\\
University of Warwick \\
Coventry, CV4 7AL England}
\email[Jones]{jdsj@maths.warwick.ac.uk}
 \date{\today}
  
\begin{document}

\begin{abstract}
     Let $M$ be a closed, oriented manifold of dimension $d$. Let $LM$ be the space of smooth loops
in $M$.   In  \cite{chsull} Chas and Sullivan defined a product on the homology  $H_*(LM)$ of degree
$-d$.  They then investigated other structure that this product induces, including a  Batalin -Vilkovisky
structure, and  a Lie algebra structure on   the $S^1$ equivariant homology
$H_*^{S^1}(LM)$.  These algebraic structures, as well as others,  came under the general heading of
the ``string topology" of $M$.  In this paper we will describe a   realization of the Chas - Sullivan loop
product  in terms of a ring spectrum structure on the Thom spectrum of a certain virtual
bundle over the loop space.  We also show that an operad action on the  homology of the loop space
discovered by  Voronov   has a homotopy theoretic realization on the level of Thom spectra.
This  is the   ``cactus operad"  defined in \cite{voronov}   
 which is equivalent to operad of framed disks in $\br^2$.  This operad action realizes the Chas - Sullivan BV
structure on  $H_*(LM)$.    We
then describe a cosimplicial model of this ring spectrum, and  by applying   the singular cochain functor to this
cosimplicial spectrum we show that this ring structure can be interpreted as the cup product in the
Hochschild cohomology,
$HH^*(C^*(M); C^*(M))$.  
  
\end{abstract}
\maketitle

\section*{Introduction  }
Let $M^d$ be a smooth, closed  $d$ - dimensional manifold, and let $LM =  C^\infty(S^1, M)$ be the
space of smooth loops in $M$.  In \cite{chsull} Chas and Sullivan described an intersection product on
the homology   $H_*(LM)$ having total degree $-d$,
$$
\circ : H_q(LM) \otimes H_r(LM) \to H_{q+r-d}(LM).
$$
In this paper we show that this product is realized by a geometric structure, not on the loop space
itself, but on the Thom spectrum of a certain bundle over $LM$. 
We describe this structure both homotopy theoretically and
simplicially, and in so doing describe the relationship of the Chas -
Sullivan product to the cup product in Hochschild cohomology.   We
now make these statements more precise.

Consider the standard parameterization of the circle by the unit interval, $exp : [0,1] \to S^1$
defined by $exp(t) = e^{2\pi i t}$.  With respect to this parameterization  we can regard a loop
$\gamma \in LM$ as a map $\gamma : [0,1] \to M$ with $\gamma (0) = \gamma (1)$.  Consider
the evaluation map 
 
\begin{align}
ev : LM &\to M \notag\\ 
\gamma &\to \gamma (1). \notag
\end{align}

Now let $\iota : M \to \br^{N+d}$ be a fixed embedding of $M$ into codimension $N$ 
Euclidean space.  Let $\nu^N \to M$ be the $N$ - dimensional normal bundle.  Let
$Th(\nu^N)$ be the Thom space of this bundle.  Recall the famous result of Atiyah
\cite{atiyah} that $Th(\nu^N)$ is  Spanier - Whitehead dual to  $M$ with a disjoint basepoint
which we denote by $M_+$.  Said more precisely, let $\mtm $ be the spectrum given by desuspending
this Thom space, 
$$
\mtm = \Sigma^{-(N+d)}Th(\nu^N).
$$
Then  there are maps of spectra
$$
S^0 \to M_+ \wedge \mtm \quad \text{and} \quad M_+ \wedge \mtm \to S^0
$$
that establish $\mtm$ as the $S - dual$ of $M_+$.  That is, these maps induce an
equivalence with the function spectrum $\mtm \simeq  Map(M_+, S^0)$.   

Notice by duality, the diagonal map $\Delta : M \to M \times M$ induces a map of spectra
$$
\Delta^* : \mtm \wedge \mtm \to \mtm
$$
that makes $\mtm$ into a ring spectrum with unit.  The unit $S^0 \to \mtm$ is the map
dual to the projection $M_+ \to S^0$.  

\med

Now let
$Th(ev^*(\nu^N))$ be the Thom space of the pull back bundle $ev^*(\nu^N) \to LM$. Define the
spectrum
$$
 \ltm = \Sigma^{-(N+d)} Th(ev^*(\nu^N)).
$$
 The  main goal of this paper is to define and study a product structure on the spectrum
$\ltm$ which  among other properties makes the evaluation map $ev : \ltm \to \mtm$
a map of ring spectra.  Here, by abuse of notation, $ev$ is referring the map of Thom
spectra induced by the evaluation map $ev : LM \to M$.  We will prove the following theorem.

\med
\begin{theorem}\label{ring}  Let $M$ be a smooth, closed manifold of dimension $d$.  The spectrum $\ltm$
is a   ring spectrum with unit, whose multiplication
$$
\mu : \ltm \wedge \ltm \la \ltm
$$
satisfies the following properties.
\begin{enumerate}
\item The evaluation map $ev : \ltm \to \mtm $ is a map of ring spectra.
\item There is a map of ring spectra $\rho: \ltm \to \Sigma^\infty (\Omega M_+)$
where the target is the suspension spectrum of the based loop space with a disjoint
basepoint.  Its ring structure is induced by the usual  product on the based loop
space.  In  the case when $M$ is oriented,   the homology of the map $\rho_*$ is given by the composition
$$
\begin{CD}
\rho_* : H_q(\ltm)  @>\tau > \cong > H_{q+d}LM  @>\iota >> H_q(\Omega M)
\end{CD}
$$
where like above, $\tau $ is the Thom isomorphism,  and the map $\iota$ takes  a $(q+d)$ -
cycle in $LM$ and intersects it with the based loop space viewed as a codimension $d$ -
submanifold. 
\item  When $M$ is orientable the  ring structure is compatible with the Chas - Sullivan homology product
in the sense that the following diagram commutes:
$$
\begin{CD}
H_q (\ltm) \otimes H_r(\ltm)  @>>> H_{q+r}(\ltm \wedge \ltm) @>\mu_* >> H_{q+r}(\ltm   )\\
@V\tau V\cong V    & & @V\cong V \tau V \\
H_{q+d}(LM) \otimes H_{r+d}(LM)  & @>> \circ >  &  H_{q+r+d}(LM)
\end{CD}
$$
\end{enumerate}
\end{theorem}

\med
\bfl 
\bf Remarks. 1.  \rm The fact that $\ltm$ is a ring spectrum was also observed by W. Dwyer and H. Miller.

\bf 2.\rm In \cite{chsull} Chas and Sullivan define a regrading of the homology  of
the loop space  
$$
\bh_q = H_{q+d}(LM)
$$
with respect to which the product $\circ$  is of total degree zero. We observe that
the Thom isomorphism defines an isomorphism
$$
\bh_* \cong H_*(\ltm)
$$ which respects gradings, and  by the above theorem is an isomorphism of rings, where the
ring structure on the right hand side comes from the ring spectrum structure of $\ltm$. 
\efl

\med
 We also show that there is a similar homotopy theoretic realization of the operad structure on the homology
of the loop space
  found by   Voronov   \cite{voronov}.   The operad involved is the  ``cactus operad"
$\cc$  which has the  homotopy type
 of  the operad of framed little disks in $\br^2$ (theorem 2.3 of \cite{voronov}).  Using a result of
Getzler \cite{getzler}, Voronov showed this operad structure  realizes the Batalin - Vilkovisky structure on
$H_*(LM)$ discovered by Chas and Sullivan.       Therefore the   following theorem can be viewed as a
homotopy theoretic realization of the $BV$ - algebra  on $H_*(LM)$. 

\med
 
\begin{theorem}\label{operad}  Let $M$ be a smooth, closed,   $d$ -dimensional manifold. Then for
each
$k$ there is a
   virtual bundle $\theta_k$over $\cc_k \times (LM)^k$ of virtual dimension $-d(k-1)$ , and a
stable map from its Thom spectrum to the loop space,
$$
\zeta_k : (\cc_k \times (LM)^k)^{\theta_k} \to LM
$$ that satisfies the following property.  If $M$ is orientable, then an orientation on $M$ induces an
orientation on the virtual bundle $\theta_k$, and  in homology, the compostion
$$
\begin{CD} 
H_{q+d(k-1)}(\cc_k \times (LM)^k) @>\tau_* >\cong > H_q((\cc_k \times (LM)^k)^{\theta_k} ) @>(\zeta_k)_*
>>H_q(LM)
\end{CD}
$$ is the   Voronov operad structure, which in turn realizes the Chas - Sullivan BV structure.  
\end{theorem}

Our final result has to do with the simplicial structure of $\ltm$, and the resulting simplicial
description of the product.

Let $S^1_*$ be the simplicial set decomposition of the circle which has one zero simplex and
one nondegenerate one simplex.  In this decomposition there are $n+1$ $n$ - simplices, all of
which are degenerate for $n > 1$.  We write this as $S^1_n = \{n+1\}$.  Now given any space
$X$, there is a resulting cosimplicial model, $ \bx_*$,  for the free loop space $LX$. 
 The   $n$ - simplices of $\bx_*$ are given by maps 
$$
\bx_n = Map (S^1_n, X) = Map(\{n+1\}, X) = X^{n+1}.
$$
The coface and codegeneracy maps of $(LX)_*$ are dual to the face and degeneracy
maps of $S^1_*$.  

Our next result states that there is a similar cosimplicial model for $\ltm$.

\med
\begin{theorem}\label{cosimp}  For $M$ a closed, oriented manifold, the spectrum $\ltm$ has
the structure of a cosimplicial spectrum which we write as $\tbx_*$.  The $k$ simplices of 
$\tbx_*$ are given by
$$
\tbx_k  = (M^k)_+ \wedge \mtm.
$$
  This cosimplicial structure
has the following properties.
\begin{enumerate}
\item  The ring structure of $\ltm$ is realized on the (co)simplicial level
by pairings
$$
\mu_k : \left( (M^k)_+ \wedge \mtm\right) \wedge \left((M^r)_+ \wedge \mtm\right) \to
(M^{k+r})_+ \wedge \mtm
$$
defined by
$$
\mu_k (x_1, \cdots , x_k; u) \wedge (y_1, \cdots , y_r; v) = (x_1, \cdots , x_k, y_1, \cdots
y_r ; \Delta^* (u \wedge v))
$$
where $\Delta^*$ is the ring structure defined on $\mtm$ described earlier.   
\item  If $M$ is orientable, then  applying the singular  chain functor $C_*( - )$ to the cosimplicial space
$\tbx_*$  yields a natural chain homotopy equivalence between the  chains of $\ltm$ and the
Hochschild cochain complex
$$
f_* : C_*(\ltm) \cong C_*(\tbx_*) \cong CH^*(C^*(M); C^*(M)).
$$
Here the notation $CH^*(A;A)$ refers to the Hochschild cochain complex 
of the form
$$
A \to Hom(A; A) \to \cdots  \to Hom(A^{\otimes n}; A) \to Hom(A^{\otimes n+1}; A) \to
\cdots
$$
Furthermore, the pairing on the chains $C_*(\ltm)$ induced by the ring spectrum structure
corresponds via the chain homotopy equivalence  $f_*$ to the cup product pairing in
$CH^*(C^*(M); C^*(M))$.  This gives   ring isomorphisms in homology,
$$
\begin{CD}
\bh_* \cong  H_*(\ltm) @>f_* >\cong >  HH^*(C^*(M); C^*(M)).
\end{CD}
$$
\end{enumerate}
\end{theorem}

\bfl
\bf Remark.  \rm The fact that the Chas - Sullivan product is realized as the cup product
in Hochschild cohomology was also observed by T. Tradler, and will
appear in his CUNY Ph.D thesis.

\efl
\med

The paper will be organized as follows.  In section 1 we will show how to realize the Chas  -
Sullivan product using the Pontrjagin - Thom constuction.  
  We will  use this to prove theorem \ref{ring}.  In section 2 we recall the definition of the cactus
operad
$\cc$ and prove theorem \ref{operad}.  In section 3 we will recall the cosimplicial study of the loop space
done by the second author in
\cite{jones}, apply the Thom spectrum construction to it, and use it to prove theorem \ref{cosimp}.  

\med
The authors are grateful to J. Klein, I. Madsen, J. Morava, G. Segal, D. Sullivan, U. Tillmann, and A. Voronov for
helpful conversations regarding this material.

\vfill \eject 
 
\section{The ring structure on $\ltm$:  The proof of theorem 1}

\med
In this section we will describe the ring spectrum structure of the Thom spectrum
$\ltm$ defined in the introduction, discuss some its properties, and prove theorem 1.
 
\med
\begin{proof}
The multiplicative structure $\mu : \ltm \wedge \ltm \to \ltm$  will be defined
using the  Pontrjagin - Thom construction.  We therefore begin by recalling some
properties of this construction.

Let
$$  e: P^k \hk  N^{n+k}$$ be an embedding of closed manifolds.
Let $\nu_e$ be a tubular neighborhood   of $e(P^k)$, which we identify with the total
space of the normal bundle of the embedding.  Let
$$
\tau : N^{n+k} \to \nu_e \cup \infty
$$
be the Pontrjagin - Thom collapse map to the one point compactification, defined by
 $$\tau (x) = \begin{cases} x \quad &\text{if}  \ x \in \nu_e  \\
\infty  \quad &\text{if} \ x \notin  \nu_e.\\
 \end{cases}$$

If we identify the  compactification with the Thom space of the normal bundle,
$\nu_e \cup \infty \cong P^{\nu_e}$, then in the oriented case we can apply the Thom
isomorphism $\begin{CD} u_*: H_q(P^{\nu_e}) @>\cong >> H_{q-n}(P^k) \end{CD}$, to
get the ``push-forward", or ``umkehr" map,
$$
\begin{CD}
e_{!} : H_q(N^k) @>\tau_* >> H_q(P^{\nu_e}) @>u_* >\cong >   H_{q-n}(P^k).
\end{CD}
$$

Recall that in the case of the diagonal embedding of a $d$ - dimensional closed  
manifold,
$$
\Delta : M^d \to M^d \times M^d
$$
that the normal bundle is isomorphic to the tangent bundle,
$$
\nu_\Delta \cong TM
$$
so that the  Pontrjagin - Thom map is a map $\tau : M \times M \to M^{TM}$. 
Furthermore, in the oriented case the push - forward map in homology,

$$
\begin{CD}
\Delta_{!} : H_*(M^d \times M^d ) @>\tau_* >> H_*(M^{TM}) @>u_* >\cong >  
H_{*-d}(M^d).
\end{CD}
$$
is simply the intersection product. 

Now the Pontrjagin - Thom construction also applies when one has a vector bundle over
the ambient manifold of an embedding.  That is, if one has an embedding
$e: P^k \hk N^{n+k}$ as above, and if one has a vector bundle (or virtual bundle)
$\zeta \to N$, then one obtains a Pontrjagin - Thom map 
$$
\tau : \zeta \cup \infty \to \nu (\iota^*(\zeta))  \cup \infty
$$
where $\nu (\iota^*(\zeta))$ is the tubular neighborhood of the induced embedding
of total spaces $\iota^*(\zeta) \hk \zeta$.  Now $\zeta \cup \infty $ is the Thom
space $N^\zeta$, and $\nu (\iota^*(\zeta))  \cup \infty$ is the Thom space
$P^{\iota^*(\zeta) \oplus \nu_e}$.  So the  Pontrjagin - Thom map  is a map 
$$
\tau : N^\zeta \to P^{\iota^*(\zeta) \oplus \nu_e}.
$$

Moreover this construction works when $\zeta$ is a virtual bundle over $N$ as well.
In this case when $\zeta  = -E$, where $E \to N$  is a $k$ - dimensional vector bundle over $N$,
then the Thom spectrum $N^\zeta = N^{-E}$ is defined as follows.  Suppose the bundle $E$ is embedded
in a $k +M$ dimensional trivial bundle, $E \hk N \times \br^{k+M}$.  Let $E^\perp$be the  $M$ -
dimensional orthogonal complement bundle to this embedding $E^\perp \to N$.  Then
$$
N^{-E} = \Sigma^{-(N+k)}N^{E^\perp}.
$$
Notice that when $E$ is oriented, the Thom isomorphism is of the form $u_* : H^q(N) \cong 
H^{q-k}(N^{-E})$.  

In particular, applying the  Pontrjagin - Thom construction    to the  diagonal embedding
$\Delta : M \hk M \times M$, using the virtual bundle $-TM \times -TM$ over $M \times M$,
we get  a map of Thom spectra,
$$
\tau : ( M \times M)^{-TM \times -TM} \to M^{TM \oplus \Delta^*(-TM \times -TM)}
$$
or,
$$
\tau : M^{-TM} \wedge M^{-TM} \to M^{-TM}.
$$
When $M$ is oriented,   this map still realizes the intersection pairing on $H_*(M)$  after applying the Thom
isomorphism.  The  map $\tau$ defines a ring spectrum structure on $\mtm$ that is well known to be the
Spanier - Whitehead dual of the diagonal map $\Delta : M \to M\times M$. 

To construct the ring spectrum  pairing $\mu : \ltm \wedge \ltm \to \ltm$, we basically ``pull back" the
structure $\tau$ over the loop space.  

To make this precise, let $ev \times ev : LM \times LM \to M \times M$ be the product of the evaluation
maps, and define $LM \times_M LM$ to be the fiber product, or pull back:

\begin{equation}\label{lmm}
\begin{CD}
\lmm  @> \tilde \Delta>> LM \times LM \\
@V ev VV  @VV ev \times ev V \\
M @>>\Delta > M \times M.
\end{CD}
\end{equation}
Notice that $\lmm$ is a codimension $d$ submanifold of the infinite dimensional manifold $LM \times
LM$,  and can be thought of as
$$
\lmm = \{ (\alpha , \beta) \in LM \times LM \ \text{such that} \ \alpha (0) = \beta (0)\}.
$$
Notice that there is also  a natural map from $\lmm$ to the loop space $LM$
defined by first applying $\alpha$ and then $\beta$.  That is,
\begin{align}\label{gamma}
\gamma : \lmm &\to LM \\
(\alpha, \beta) & \to \alpha * \beta \notag
\end{align}
where 
$$
\alpha *\beta (t) =  \begin{cases}  \alpha(2t) \quad &\text{if} \ 0\leq t \leq \frac{1}{2} \\
\beta(2t-1) \quad &\text{if} \ \frac{1}{2} \leq t \leq 1.\\
 \end{cases}
$$
Notice that when restricted to the product of the based loop spaces,
$\Omega M \times \Omega M \subset \lmm$, then $\gamma$ is just the $H$ - space product on the
based loop space, $\Omega M \times \Omega M \to \Omega M$. 

\bfl

\bf Remark.  \rm  This definition needs to be modified slightly since $\alpha*\beta$ may not be smooth
at $1/2$ or $1$.  This is dealt with in a standard way by  first ``dampening"  $\alpha$ and $\beta$ by 
using a smooth bijection
$\phi : [0,1] \to [0,1]$ with the property that all of its derivatives approach zero at the endpoints $\{0\}$
and
$\{1\}$, to reparameterize $\alpha$ and $\beta$.  This will allow them to be ``spliced" together by the
above formula without losing any smoothness.  This is a standard construction, so we leave its details to the
reader.  

\efl

Notice that by its definition (\ref{lmm}) the  embedding $\tilde \Delta : LM \times_M LM \hk LM \times
LM$  has a tubular neighborhood $\nu (\tilde \Delta)$ defined to be the inverse image of the tubular
neighborhood of the diagonal $\Delta : M \hk M \times M$:
$$
  \nu (\tilde \Delta) = ev^{-1}(\nu(\Delta)).
$$
Therefore this neighborhood is homeomorphic to the total space of the $d$  - dimensional vector bundle
given by pulling back the normal bundle of the embedding $\Delta$, which is the tangent bundle of $M$:
$$
\nu (\tilde \Delta) \cong ev^*(\nu_\Delta) = ev^*(TM).
$$
Thus there is a Pontrjagin - Thom collapse map
$$
\tau : LM \times LM \to \lmm^{ev^*(TM)}.
$$
As described earlier, we ease the notation by refering to this Thom spectrum as $(\lmm)^{TM}$.
By the naturality of the Pontrjagin - Thom construction, we have a commutative diagram of spectra,

\begin{equation}\label{specev}
\begin{CD}
LM \times LM   @>\tau >> (\lmm)^{TM} \\
@Vev VV   @VV ev V \\
M \times M @>>\tau >  M^{TM}
\end{CD}
\end{equation}

 Since in the oriented case $\tau_* : H_*(M \times M) \to H_*(M^{TM}) \cong H_{*-d}(M)$ is the
intersection product, then  
$$
\begin{CD}
H_*(LM \times LM) @>\tau_* >> H_*((\lmm)^{TM}) @>u_* >> H_{*-d}(\lmm)
\end{CD}
$$
 can be viewed (as is done in Chas  - Sullivan \cite{chsull}) as taking a cycle in $LM \times LM$, and 
``intersecting"  it with the codimension $d$ submanifold $\lmm$. 

Now observe that the map $\gamma : \lmm \to LM$ defined above (\ref{gamma})  preserves the
evaluation map.  That is, the following diagram commutes:
\begin{equation}\label{ev}
\begin{CD}
\lmm @>\gamma >>  LM \\
@Vev VV  @VV ev V \\
M  @>> = > M
\end{CD}
\end{equation}
Thus $\gamma$ induces a map of bundles $\gamma : ev^*(TM) \to ev^*(TM)$, and therefore a map of
Thom spectra,
$$
\gamma : (\lmm)^{TM} \to LM^{TM}.
$$
Now consider the composition
\begin{equation}\label{tmu}
\begin{CD}
\tilde \mu : LM \times LM  @>\tau >> (\lmm )^{TM} @>\gamma >> LM^{TM}
\end{CD}
\end{equation}
In the oriented case the homomorphism
\begin{equation}\label{prod}
\begin{CD}
H_*(LM \times LM)  @>\tilde \mu_* >> H_{*}(LM^{TM}) @>u_* > \cong >  H_{*-d}(LM)
\end{CD}
\end{equation}
takes a cycle in $LM \times LM$, intersects in with the codimension $d$ - submanifold $\lmm$, 
and maps it via $\gamma$ to $LM$.  This is the definition of the Chas - Sullivan product 
 $H_*(LM)$.

Now as we did before with the diagonal embedding, we can perform the Pontrjagin - Thom construction
when we pull back the virtual bundle $-TM \times -TM$ over $LM \times LM$.  That is, we get a map of
Thom spectra
$$
\tau : (LM \times LM)^{(ev \times ev)^*(-TM \times -TM)} \la  (\lmm)^{ev^*(TM) \oplus
ev^*(\Delta^*(-TM \times -TM))}.
$$
But since $ev^*(\Delta^*(-TM \times -TM)) = ev^*(-2TM)$, we have
$$
\tau : \ltm \wedge \ltm \la (\lmm)^{TM \oplus -2TM} = (\lmm)^{-TM}.
$$
Now by the  commutativity of (\ref{ev}), $\gamma$ induces a map of Thom spectra,
$$
\gamma : (\lmm)^{-TM} \to \ltm
$$
and so we can define the ring structure on the Thom spectrum $\ltm$ to be the composition
 \begin{equation}\label{mu}
\begin{CD}
\mu : \ltm \wedge \ltm @> \tau >> (\lmm )^{-TM} @>\gamma >> \ltm.
\end{CD}
\end{equation}

A few properties of this map $\mu$ are now immediately verifiable.

First, $\mu$ is associative.  This follows from the naturality of the
Pontrjagin - Thom construction, and the fact that the map $\gamma$
is associative.  (Strictly speaking, formula (\ref{gamma}) is $A_\infty$ -
associative as is the usual formula for the product on the based loop
space, $\Omega M$.  However the standard trick of replacing $\Omega
M$ with ``Moore loops" changes the $A_\infty$ structure to a strictly
associative structure.  The same technique applies to the map
$\gamma$.  Otherwise, the spectrum $\ltm$  will have the structure of
an $A_\infty$ ring spectrum.)

Also, notice tha $\ltm$ has a unit, $\iota : S^0 \to \ltm$, defined by
the composition 
$$
\begin{CD}
\iota : S^0 @> j >> \mtm @>\sigma >> \ltm
\end{CD}
$$
where $j $    is the unit of the ring spectrum structure of
$\mtm$, and $\sigma$ is  the map of Thom spectra induced by the
section of the evaluation map $ev : LM \to M$ defined by viewing
points in $M$ as constant loops.

Notice furthermore that in the oriented case, after applying the Thom
isomorphism, $\mu_*$ induces the same homomorphism as $\tilde
\mu_*$, and so by  (\ref{prod}) the following diagram commutes:
\begin{equation}
\begin{CD}
H_{q-2d}(\ltm \wedge \ltm) @>\mu_* >> H_{q-2d}(\ltm) \\
@Vu_* V \cong V    @V\cong V u_*V \\
H_q(LM \times LM)  @>>\circ >  H_{q-d}(LM)
\end{CD}
\end{equation}
where $\circ : H_q(LM \times LM)  \to  H_{q-d}(LM)$ is the
Chas - Sullivan product.   This proves part (3) of theorem
\ref{ring}.

Now by the naturality of the Pontrjagin - Thom construction, the
following diagram of Thom spectra commutes (compare (\ref{specev}))

\begin{equation}
\begin{CD}
\ltm \wedge \ltm @>\mu >>  \ltm \\
@Vev \times evVV   @VV ev V \\
\mtm \wedge \mtm @>>\tau > \mtm.
\end{CD}
\end{equation}
Thus the evaluation map $ev : \ltm \to \mtm$ is a map of ring spectra,  which proves part
1 of theorem  \ref{ring}.

\med
We now verify part 2 of theorem \ref{ring}.  Let $x_0 \in M$ be a base point, and consider
the following pullback diagram:
\begin{equation}\label{omega}
\begin{CD}
\Omega M  @>j >> LM \\
@VpVV  @VVev V \\
x_0 @>> i > M.
\end{CD}
\end{equation}
Thus the embedding $j : \Omega M \hk LM$ is an embedding of a codimension
$d$ submanifold, and has tubular neighborhood, $\nu (j)$ equal to the inverse image of the
tubular neighborhood of the inclusion of the basepoint $i : x_0 \hk M$.  This tubular
neighborhood is simply a disk $D^d$, and so 
$$
\nu (j) \cong \Omega M \times D^d  \cong \epsilon^d
$$
where $\epsilon^d$ reflects the $d$ dimensional trivial bundle over $\Omega M$.  Thus the
Thom - Pontrjagin construction makes sense, and is a map
$$
\tau : LM \to \Omega M^{\epsilon^d} = \Sigma^d(\Omega M_+)
$$
where the last space is the $d$ - fold suspension of $\Omega M$ with a disjoint
basepoint.  In homology, the homomorphism
$$
\tau_* : H_q(LM) \to H_{q}(\Sigma^d(\Omega M_+) = H_{q-d}(\Omega M)
$$
denotes the map that is obtained by intersecting a $q$ -cycle in $LM$ with the codimension
$d$ submanifold $\Omega M$.  

By performing the Pontrjagin - Thom construction after pulling back the virtual bundle $-TM$
over $LM$, we get a map of Thom spectra
$$
\tau : \ltm \to (\Omega M)^{\epsilon^d \oplus j^*ev^*(-TM)}.
$$
But by the commutativity of diagram (\ref{omega}), $j^*ev^*(-TM) = p^*i^*(-TM)$, which is
the trivial, virtual  $-d$ dimensional bundle which we denote $\epsilon^{-d}$.  So the 
Pontrjagin - Thom map is therefore a map of spectra
\begin{equation}\label{rho}
\begin{CD}
\rho : \ltm @>\tau >> (\Omega M)^{\epsilon^d \oplus j^*ev^*(-TM)} = (\Omega
M)^{\epsilon^d \oplus \epsilon^{-d}} = \Sigma^\infty (\Omega M_+)
\end{CD}
\end{equation}
where  $\Sigma^\infty (\Omega M_+)$ denotes the suspension spectrum of
the based loop space of $M$ with a disjoint basepoint.   To complete the proof of theorem
\ref{ring} we need to prove that $\rho : \ltm \to \Sigma^\infty (\Omega M_+)$ is a map of
ring spectra.  Toward this end, consider the following diagram of pull back squares:

\begin{equation}
\begin{CD}
\Omega M \times \Omega M @>\iota >> LM \times_M LM @>\tilde \Delta >> LM \times LM \\
@VVV   @VVev V  @VVev \times ev V \\
x_0 @>>i > M @>> \Delta >  M \times M.
\end{CD}
\end{equation}
This gives Pontrjagin - Thom maps
\begin{equation}\label{tautau}
\begin{CD}
\ltm \wedge \ltm  @>\tau >> (\lmm)^{-TM} @> \tau >> \Sigma^\infty (\Omega M_+) \wedge
\Sigma^\infty (\Omega M_+). 
\end{CD}
\end{equation}
Notice that by the naturality of the Pontrjagin - Thom construction, the above composition is
equal to
$$
\rho \wedge \rho : \ltm \wedge \ltm \la \Sigma^\infty (\Omega M_+) \wedge
\Sigma^\infty (\Omega M_+). 
$$

Now notice that by the formula for the map $\gamma : \lmm \to LM$, the following
diagram commutes:
$$
\begin{CD}
\Omega M \times \Omega M  @>\iota >>  \lmm \\
@V m VV  @VV \gamma V \\
\Omega M  @>> \subset > LM
\end{CD}
$$
where $m$ is the usual multiplication on the based loop space. 
Pulling back the virtual bundle $-TM$ over $LM$, and applying the Pontrjagin - Thom
construction, we then get a commutative diagram of spectra,
\begin{equation}
\begin{CD}
\ltm \wedge \ltm @>\tau >> (\lmm )^{-TM}  @>\tau >> \Sigma^\infty (\Omega M_+) \wedge
\Sigma^\infty (\Omega M_+)\\
& &@V \gamma VV  @VV m V \\
& & \ltm @>>\rho >  \Sigma^\infty (\Omega M_+) 
\end{CD}
\end{equation} 
Now as observed above (\ref{tautau}),  the top horizontal composition $\tau \circ \tau$
is equal to $\rho \wedge \rho : \ltm \wedge \ltm \to \Sigma^\infty (\Omega M_+) \wedge
\Sigma^\infty (\Omega M_+)$.   Also, $\gamma \circ \tau$ is, by definition, the ring
structure  $\mu : \ltm \wedge \ltm \to \ltm$.  Thus the following diagram of spectra
commutes:
$$
\begin{CD}
\ltm \wedge \ltm @>\rho \wedge \rho >> \Sigma^\infty (\Omega M_+) \wedge
\Sigma^\infty (\Omega M_+)\\
@V\mu VV   @VV m V \\
\ltm @>> \rho >  \Sigma^\infty (\Omega M_+).
\end{CD}
$$
Thus $\rho$ is a map of ring spectra, which completes the proof of theorem \ref{ring}.
\end{proof}

\vfill \eject

\section{The operad structure}

In this section we describe a homotopy theoretic realization of the operad structure on the homology of the
loop space
  found  Voronov \cite{voronov}.  
The operad involved is the  ``cactus operad"
$\cc$ which has the  homotopy
 of  the operad of framed little disks in $\br^2$ (theorem 2.3 of \cite{voronov}).  According to Getzler's
result
\cite{getzler} this is precisely what is needed to induce the Batalin - Vilkovisky algebra
structure, and Voronov observed that his operad structure induced the BV structure on $H_*(LM)$ described
by Chas and Sullivan. Therefore
 the constructions in this section  can be viewed as a homotopy theoretic realization of  Chas and Sullivan's
$BV$ - algebra structure  on $H_*(LM)$.  

We begin by recalling the definition of   cactus operad $\cc$.  
 A point in the  space $\cc_k$  is a
collection of $k$ oriented, parameterized circles $c_1,
\cdots , c_k$,  with radii $r_i$ so that $\sum_{i-1}^k r_i = 1$. Each circle has a marked point
$x_i \in c_i$ given by the image under the parameterization of the basepoint $1 \in S^1$.  Moreover the
circles can intersect each other at a finite number of points  (vertices)  to create a ``cactus - type
configuration".  Strictly speaking this means that the dual graph of this configuration is a tree.  That is, the
``cactus" (i.e the union of the circles) must be connected and have no ``extra loops".   (This is the tree
condition on the dual graph.) The boundary of the cactus (that is the union of the circles) is also equipped
with a fixed basepoint $y_0$ together with a choice of which component the basepoint $y_0$ lies in. Say
$y_0 \in c_{j_0}$. 
(This choice is only relevant if the basepoint happens to be one of intersection points.) The
edges coming into any vertex are also equipped with a cyclic ordering. An ordering of the circle components
making up $c$ as well as the vertices within each component is part of the data making up the cactus.  The
topology of the space of cacti with $k$ - components, $\cc_k$ is described in \cite{voronov}.

 \med
 \begin{center}
  \scalebox{.50}{\includegraphics{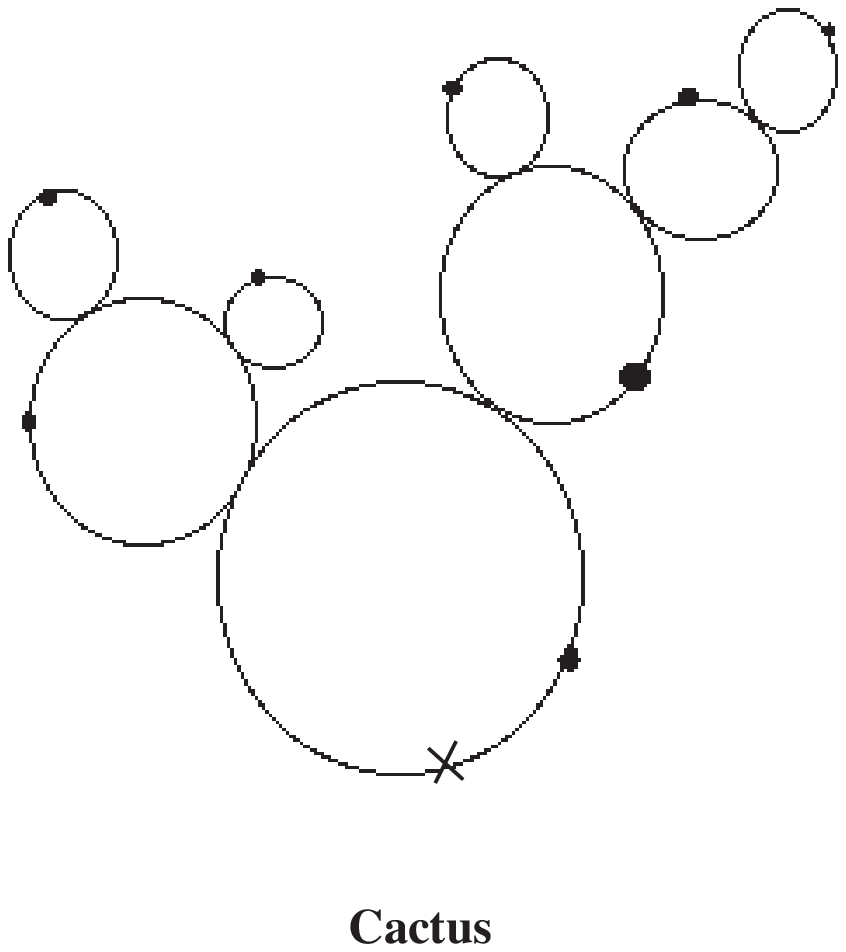}} 
   \end{center}

Notice that a cactus  (i.e a point in $\cc_k$) comes equipped with a well defined map from
the unit circle to the boundary of the cactus.  That is, the map begins at the  basepoint $y_0 \in c_{j_0}$, 
  then traverses the circle $c_{j_0}$ in the direction of its orientation in a
length preserving manner.  When a point of intersection with another circle is reached,
the loop then traverses that circle, in the direction of its orientation. This path is
continued until it eventually arrives back at the original basepoint $y_0 \in c_{j_0}$.  Given a
cactus $c = \langle c_1, \cdots , c_k \rangle \in \cc_k$ we let $\delta_c : S^1 \to c_1
\cup \cdots \cup c_k$ be this loop that traverses the boundary of the cactus.

As mentioned above, the result that the  cactus operad is homotopy equivalent to the  framed little disk
operad is a theorem of Voronov \cite{voronov}.  Roughly, the relationship between these operads can seen as
follows. 

Let $\cf_k$ be the $k^{th}$ space of $k$ framed little disks in $\br^2$.  Recall that there is a homotopy
equivalence,
$$
\cf_k \simeq F(\br^2, k) \times (S^1)^k 
$$
where $F(\br^2,k)$ is the configuration space of $k$ ordered, distinct points in $\br^2$. 
One can construct a map $\phi_k : \cc_k \to \cf_k$ in the following way.  View a cactus
$c \in \cc_k$ as a planar graph.       From this perspective, the center
of each component circle determines a configuration of $k$ points in the plane, $\phi^1_k (c) \in F(\br^2,k)$.
Furthermore, the parameterizations of the component circles determine a point $\phi^2_k(c) \in (S^1)^k$.
Then $\phi_k$ is defined (up to homotopy) as the composition
$$
\begin{CD}
\phi_k: \cc_k @>\phi^1_k \times \phi^2_k >> F(\br^2, k) \times (S^1)^k @>\simeq >> \cf_k.
\end{CD}
$$
One can see that the maps $\phi_k$ are homotopy equivalences by an inductive argument on the number of
components $k$, using natural fibrations $\cc_k \to \cc_{k-1}$ defined by ignoring the last component
circle that does not contain the basepoint, and  $\cf_k \to \cf_{k-1}$ defined by ignoring the last disk. 
We refer the reader to Voronov's work for a complete proof of this homotopy equivalence.

Notice that $\cc_k$ has a free action of the symmetric group $\Sigma_k$ defined by
permuting the ordering of the circles.  The operad action
$$
\xi : \cc_k \times \left(\cc_{j_1} \times \cdots \times \cc_{j_k}\right) \to \cc_j
$$
where $j = j_1 + \cdots + j_k$ is defined as follows.  Let $c \times (c^1, \cdots , c^k) \in
\cc_k \times \left(\cc_{j_1} \times \cdots \times \cc_{j_k}\right)$.  Scale the cactus
$c^1$ down so that its total radius is the radius $r_1$ of the first component $c_1$ of
$c \in C_k$.   Similarly, scale each of the cacti $c^i$ down so that its total radius is the
radius $r_i$ of the $i^{th}$ component $c_i$ of $c$.  By using the loops 
$\delta_{c^i}$ (scaled down appropriately) we identify the component circles $c_i$ with
the boundary of the cactis $c^i$.  This has the effect of replacing the $k$ component
circles $c_1, \cdots , c_k$ making up the cactus $c$, by the $k$ cacti, $c^1, \cdots ,
c^k$.  This produces a cactus with $j_1 + \cdots +j_k = j$ component circles. 

\med
Our goal  in this section is to prove theorem \ref{operad}  as stated in the introduction.

\begin{proof}  Given a cactus $c = \langle c_1, \cdots , c_k\rangle \in \cc_k$,  
define $L_cM$ to be the mapping space  
$$ L_cM = Map (c, M).$$ This space consists of maps from the union $c_1 \cup \cdots \cup
c_k
\to M$.  The map from the circle $\delta_c : S^1 \to c_1 \cup \cdots \cup c_k$ defines a
map from $L_cM $ to the loop space,
\begin{align}\label{gammac}
\gamma_c : L_cM &\to LM \\  
f &\to f\circ \delta_c. \notag
\end{align}

Now $L_cM$ can also be viewed as the pullback of  an evaluation
mapping of the product $(LM)^k$ defined as follows.   For each component of the cactus
$c_i$, let $y_1, \cdots , y_{m_i}$ denote the vertices of $c$ that lie on $c_i$.
   $m_i =
m_i (c)$ is the number of these vertices.     Let $m_c = m_1 + \cdots + m_k$.  We define an
evaluation map
$$
ev_c : (LM)^k \la (M)^{m_c}
$$
as follows.  Let
$s_i : S^1 \to c_i$ be the identification of the unit circle with $c_i$ obtained by
scaling down the unit circle so as to have radius $r_i$, and rotating it so the basepoint
$1 \in S^1$ is mapped to the marked point $x_i \in c_i$.  Let $u_1, \cdots , u_{m_i}$ be
the points on the unit circle corresponding to $y_1, \cdots y_{m_i} \in c_i$ under the map
$s_i$.  Define 
\begin{align}
ev_{c_i} : LM  &\to (M)^{m_i} \\
\sigma &\to (\sigma (u_1), \cdots , \sigma (u_{m_i})) \notag
\end{align}
Now define
\begin{equation}\label{evc}
ev_c = ev_{c_1} \times \cdots \times ev_{c_{m_i}}  : LM \to (M)^{m_1} \times \cdots
\times (M)^{m_i} = (M)^{m_c}
\end{equation}

Now let $w_1, \cdots, w_{n_c} \in c_1 \cup \cdots \cup c_k$  be the ordered collection of vertices in the
cactus $c$ (ordered by the ordering of the components $c_1, \cdots c_k$ and the ordering of the
vertices in each component as described above).   For each
such vertex
$w_i$, let
$\mu_i$ be the number of components of the cactus on which $w_i$ lies.  We think of $\mu_i$ as the  
``multiplicity" of the vertex $w_i$.   Notice  that we have the relation
\begin{equation}
\sum_{i = 1}^{n_c} \mu_i = m_c.
\end{equation}
 The ``tree" condition on the dual of the cactus also imposes the following relation:
\begin{equation}\label{mcnc}
m_c - n_c = k - 1.
\end{equation}

Now consider the diagonal mapping 
$$
\Delta_c :  (M)^{n_c} \la (M)^{m_c}
$$
defined by the composition
$$
\begin{CD}
\Delta_c : (M)^{n_c}  @>\Delta_{\mu_1} \times \cdots \times \Delta_{\mu_{n_c}} >>
(M)^{\mu_1}\times \cdots \times (M)^{\mu_{n_c}} @> = >>  (M)^{m_c}
\end{CD}
$$
where $\Delta_{\mu_i} : M \to  (M)^{\mu_i}$ is the $\mu_i$  -fold diagonal.  Observe
that the following is a cartesian pull - back square:

\begin{equation}
\begin{CD}
L_cM  @>\tilde \Delta_c>>   (LM)^k \\
@Vev_{int} VV   @VVev_c V \\
(M)^{n_c} @>>\Delta_c >  (M)^{m_c}
\end{CD}
\end{equation}
where $ev_{int} : L_cM \to (M)^{n_c}$ evaluates a map $f : c \to M$ at the $n_c$
vertices $w_1, \cdots , w_{n_c} \in c$. 

The normal bundle $\nu (\Delta_c)$ of the diagonal embedding  
$$
\begin{CD}
\Delta_c  :   (M )^{n_c} @>\Delta_{\mu_1} \times \cdots \times \Delta_{\mu_{n_c}} >>
(M)^{\mu_1}  \times \cdots \times  (M)^{\mu_{n_c}}@> = >>   (M)^{m_c} 
\end{CD}
$$
is equal to 
$$
(\mu_1 -1) TM \times \cdots \times (\mu_{n_c}-1)TM \la (M)^{n_c},
$$
where $(q)TM$ denotes the  $q$ - fold direct sum of $TM$ with itself as a bundle over $M$.
By the tubular neighborhood theorem, we have therefore proven the following.

\begin{lemma}\label{nbd}The image of the embedding $\tilde \Delta_c : L_c M \hk (LM)^k$ has an open
neighborhood homeomorphic to the total space of the pullback $ev_{int}^*((\mu_1 -1) TM \times \cdots \times
(\mu_{n_c}-1)TM )  $.
\end{lemma}

\med
  We now consider these constructions in a parameterized way, by letting $c \in \cc_k$ vary.   Namely,
 let 
\begin{equation} \label{lkm}
L_kM = \{(c,f) : c \in \cc_k \ \text{and}  \ f \in L_cM\} 
\end{equation}
We then have a map
\begin{align}
\tilde \Delta_k : L_kM  &\to \cc_k \times (LM)^k  \notag\\
(c,f) &\to (c, \tilde \Delta_c (f)). \notag
\end{align}  We also have a map
\begin{align}
\gamma_k : L_kM  &\to  LM   \notag\\
(c,f) &\to \gamma_c (f), \notag
\end{align} where $\gamma_c$ is as defined in (\ref{gammac}).

Lemma \ref{nbd} allows us to identify an open neighborhood of the image of $\tilde \Delta_k$, as follows.
Let $\xi_k$ be the $d(k-1)$ dimensional vector bundle over $L_kM$ whose fiber over $(c,f)$ is the sum of
tangent spaces, 
$$
(\xi_k)_{|_{(c,f)}} = \bigoplus_{j=1}^{n_c} (\mu_j -1)T_{f(w_j)}M
$$
where, as above, $w_1, \cdots, w_{n_c}$ are the vertices of $c$, and $\mu_j = \mu_j(c)$is the multiplicity
of the vertex $w_j$.   It is immediate that $\xi_k$ is a well defined vector bundle.  Lemma \ref{nbd} now
yields the following.

\med
\begin{lemma}\label{bignbd}  The image of the embedding $\tilde \Delta_k : L_kM \to \cc_k \times (LM)^k$
has an open neighborhood homeomorphic to the total space of $\xi_k$. 
\end{lemma}

\med
Notice that this will allow us to perform a Pontrjagin - Thom collapse map
$$
\tau : \cc_k \times (LM)^k  \la (L_kM)^{\xi_k}.
$$
But in order to prove theorem \ref{operad}, we will need to twist this map by a virtual bundle  $\theta_k$
over $\cc_k \times (LM)^k$.   To define this bundle, first fix an embedding
$M \hk \br^{d+L}$ with $L >> d$, having normal bundle $\eta$ of dimension $L$.  Define a vector bundle
$\tilde \theta_k$ of dimension $L(k-1)$ over $\cc_k \times (LM)^k$ whose fiber at $ (c; f_1, \cdots f_k) $
is given by
$$
(\tilde \theta_k)_{|_{ (c; f_1, \cdots f_k) }} = \bigoplus_{j=1}^{n_c} (\mu_j -1) \eta_{f(x_{j,1})}
$$
where $x_{j,1}, \cdots , x_{j, \mu_j}$ are the   points on $\coprod_k S^1 $ that
map to $u_j$ under the parameterization $\begin{CD} S^1 \sqcup \cdots \sqcup S^1 @>s_1 \sqcup \cdots
\sqcup s_k >> c_1 \cup \cdots \cup c_k \end{CD}$. (Notice that by the ``tree" condition on the cactus,
there is at most one of the $x_{j,r}'s$ on any one circle in $\coprod_k S^1 $.)  

It is clear that these fibers fit together to give a vector bundle $\tilde \theta_k$ over $\cc_k \times
(LM)^k$.   Lemma \ref{bignbd} implies the following result.

\begin{lemma}\label{totnbd}  Consider the embedding of total spaces of bundles
$$
\begin{CD}
\tilde \Delta_k^*(\tilde \theta_k)  @>\tilde \Delta_k >\hk >  \tilde \theta_k \\
@VVV    @VVV \\
L_k M   @>\hk >\tilde \Delta_k >   \cc_k \times (LM)^k.
\end{CD}
$$
then the image of $\tilde \Delta_k^*(\tilde \theta_k)$ in $\tilde \theta_k$ has an open neighborhood
homeomorphic to the total space of the bundle $\tilde \Delta_k^*(\tilde \theta_k)  \oplus \xi_k$ over
$L_kM$.
\end{lemma}

This implies we have a Pontrjagin - Thom collapse map, $$\tau_k : (\cc_k \times (LM)^k)^{\tilde \theta_k}
\la (L_kM)^{\tilde \Delta_k^*(\tilde \theta_k)  \oplus \xi_k}.$$

Now by definition,  $\tilde \Delta_k^*(\tilde \theta_k)  \oplus \xi_k$ is a vector bundle of dimension
$(d+L)(k-1)$ whose fiber at $(c,f)$ is the sum of the vector spaces,
$$
\bigoplus_{j=1}^{n_c} (\mu_j -1)\left( T_{f(w_j)}M \oplus \eta_f(w_j)\right).
$$
But since $TM \oplus \eta$ is canonically trivial, we have an induced trivialization
$$
\tilde \Delta_k^*(\tilde \theta_k)  \oplus \xi_k  \cong L_kM \times \br^{(d+L)(k-1)}.
$$
Therefore the Pontrjagin - Thom collapse  is a map 
$$
\tau_k : (\cc_k \times (LM)^k)^{\tilde \theta_k}
\la  \Sigma^{(d+L)(k-1)}(L_kM_+).
$$

Let $\theta_k$ be the virtual bundle  over $L_kM$, 
$$
\theta_k = \tilde \theta_k - [(d+L)(k-1)]
$$ where $[q]$ refers to the $q$ - dimensional trivial bundle.   Notice that $\theta_k$ has virtual
dimension $-d(k-1)$, and its Thom spectrum is the $-(d+L)(k-1)$ - fold desuspension of the Thom spectrum of
$\tilde \theta_k$.  We then define $\mu_k$ to be the  composition,
\begin{equation}\label{mook}
\begin{CD}
\mu_k: (\cc_k \times (LM)^k)^{\theta_k} @> \Sigma^{-(d+L)(k-1)} \tau_k >> L_kM @> \gamma  _k >> LM.
\end{CD}
\end{equation}

Now recall from (\cite{voronov}) that the cactus operad structure on $H_*(LM)$ is given by the composition
$$
\begin{CD}
H_*(\cc_k \times (LM)^k) @> (\tilde \Delta_k)_! >> H_*(L_kM) @>(\gamma_k)_*>> H_*(LM)
\end{CD}
$$
where $(\tilde\Delta_k)_!$ is the push -forward map. (Note that our notation is different than that used
in \cite{voronov}.)    The theorem now follows from the well known fact that the Pontrjagin - Thom collapse
map realizes the push -forward
 map in homology.  \end{proof}

\vfill \eject
 
\section{A cosimplicial description of $LM$ and $\ltm$ and a proof of
Theorem \ref{cosimp}}

\med 
In this section we describe a cosimplicial  model for the spectrum  $\ltm$.    We then describe
the ring spectrum structure simplicially.  This cosimplicial model will then give a natural way of
relating the  singular  chains $C_*(\ltm) $ to the Hochschild cochain complex $CH^*(C^*(M),
C^*(M))$, and in particular relate the simplicial model for the ring structure of $\ltm$ to the cup
product structure in this cochain complex.  This will allow us to prove theorem \ref{cosimp}.  

We begin by reviewing the cosimplicial model of the loop space $LX$ for any space $X$, coming
from  a simplicial decomposition of  the circle $S^1$.  We refer the reader to \cite{jones}
for details.

Let $S^1_*$ be the simplicial set decomposition of the circle which has one zero simplex and
one nondegenerate one simplex.  In this decomposition there are $n+1$ $n$ - simplices, all of
which are degenerate for $n > 1$.  We write this as $S^1_n = \{n+1\}$.  Now given any space
$X$, there is a resulting cosimplicial model for the free loop space, $LX $, which we call
$\bx_*$.  The   $n$ - simplices of $\bx*$ are given by maps 
$$
\bx_n = Map (S^1_n, X) = Map(\{n+1\}, X) = X^{n+1}.
$$
The coface and codegeneracy maps of $\bx_*$ are dual to the face and degeneracy
maps of $S^1_*$.  They are given by the formulas

\begin{align}\label{coface}
\delta_i(x_0, \cdots , x_{n-1})  &= (x_0, \cdots, x_{i-1}, x_i, x_i, x_{i+1}, \cdots , x_{n-1}),
\quad  0 \leq i \leq n-1 \\
\delta_n(x_0, \cdots , x_{n-1}) &= (x_0, x_1, \cdots , x_{n-1}, x_0) \notag\\
\sigma_i(x_0, \cdots , x_{n+1}) &= (x_0, \cdots, x_i, x_{i+2}, \cdots, x_{n+1}), \quad 0 \leq i \leq n
\notag 
\end{align}

  Since the
geometric realization of
$S^1_*$ is homeomorphic to the circle,
$$
S^1 \cong |S^1_*|
,$$
the ``total complex"  or geometric corealization  of $\bx_*$ is homeomorphic to the
loop space,
$$
LX \cong Tot(\bx_*).
$$
This was studied in detail by the second author in \cite {jones}, and in particular the following
interpretation of this result was given.  For each $k$, let $\Delta^k$ be the standard $k$ -
simplex:
$$
\Delta^k = \{(x_1, \cdots , x_k) : 0 \leq x_1  \leq x_2  \leq  \cdots  \leq x_k \leq 1\}.
$$
Consider the maps
\begin{align}\label{fk}
f_k : \Delta^k \times LX &\la X^{k+1}   \\
(x_1, \cdots , x_k) \times \gamma &\to (\gamma (0), \gamma(x_1), \cdots , \gamma(x_k)
 ).
\notag
\end{align}

Let $\bar f_k : LX \to Map(\Delta^k, X^{k+1})$ be the adjoint of $f_k$.  Then the following
was proven in \cite{jones}.

\med
\begin{theorem}\label{jones}
Let $X$ be any space, and let $f : LX \la \prod_{k \geq 0} Map (\Delta^k, X^{k+1})$ be the
product of the maps
$\bar f_k$.  Then $f$ is a homeomorphism onto its image.  Furthermore, the image consists of
sequences of maps $\{\phi_k\}$ which commute with the coface and codegeneracy operators, which is the
total space, $Tot (\bx_*)$. 
\end{theorem}

\med
 By applying singular cochains to the maps $f_k$, one obtains maps 
$$
f_k^* : C^*(X)^{\otimes k+1} \to C^{*-k}(LX).
$$
The following was also observed in \cite{jones}.

\med
\begin{theorem}\label{hoch}
For any space $X$, the homomorphisms $f_k^* : C^*(X)^{\otimes k+1} \to C^{*-k}(LX)$ fit
together to define
  a chain map from the Hochschild complex
of the cochains of $X$ to the cochains of the free loop space,
$$
f^*: CH_*(C^*(X)) \to C^*(L(X))
$$
which is a chain homotopy equivalence when $X$ is simply connected.  Hence it    induces
an isomorphism in homology
$$
\begin{CD}
f^* : HH_*(C^*(X))  @>\cong >> H^*(L(X)).
\end{CD}
$$
\end{theorem}

\med
\bfl
\bf Remark.  \rm Let us clarify some notation.  Given an algebra (or differential graded
algebra) $A$, the the Hochschild complex of $A$,  $CH_*(A)$  is a complex of the  form
$$
\begin{CD}
\cdots @>b >> A^{\otimes n+2} @> b >> A^{\otimes n+1} @>b >> \cdots @>b>> A\otimes A @>b>>
A.
\end{CD}
$$
The homology of this algebra is denoted $HH_*(A)$. More generally if $M$ is a bimodule
over $A$, we denote by $CH_*(A; M)$ the Hochschild complex of the form
$$
\begin{CD}
\cdots @>b >> A^{\otimes n+1}\otimes M @>b> > A^{\otimes n}\otimes M @>b >> \cdots @>b>>
A\otimes M @>b>> M.
\end{CD}
$$
The homology of this complex is denoted $HH_*(A;M)$.  So in particular if $M =A$ we see
that $HH_*(A; A) = HH_*(A)$.  Dually, we denote by $CH^*(A; M)$ the Hochschild cochain
complex of the form
$$
\begin{CD}
M @>b^* >> Hom(A; M) @>b^*>>  \cdots @>b^*>> Hom(A^{\otimes n}; M) @>b^* >>
Hom(A^{\otimes n+1}; M) @>b^* >> \cdots
\end{CD}
$$
Its cohomology is denoted $HH^*(A;M)$.  By dualizing theorem \ref{jones} we obtain the
following.

\efl

\med
\begin{corollary} \label{chain}For any simply connected space $X$, there is a chain homotopy
equivalence from the singular chains of the loop space to the Hochschild cochain complex
$$
f_* : C_*(LX) \to CH^*(C^*(X); C_*(X))
$$
and so an isomorphism in homology,
$$
\begin{CD}
f_*: H_*(LX) @>\cong >> HH^*(C^*(X); C_*(X)).
\end{CD}
$$
\end{corollary}

\med
Notice that the cochain complex $CH^*(C^*(X); C_*(X))$ does \sl not \rm in general have a
natural product structure.  This is because the coefficients, $C_*(X)$, is not in general a ring.
Notice however that the Hochschild complex $CH^*(C^*(X), C^*(X))$ does in fact have a cup
product coming from the algebra structure of $C^*(X)$.  Of course when $X$ is a closed,
oriented manifold of dimension $d$,  Poincare duality gives a chain homotopy equivalence, $C_*(X)
\cong C^{d-*}(X)$, and so the cochain complex $CH^*(C^*(X); C_*(X))$ inherits an algebra
structure.  Therefore by the above corollary, $H_*(LX)$ inherits an algebra structure in this case.
We will see that this indeed realizes the Chas - Sullivan product.  We will show this by
showing that  when $M$ is a closed    $d$ - manifold, the Thom spectrum $\ltm$ inherits
a cosimplicial structure   for which the analogue of theorem \ref{hoch}  will yield
a natural chain homotopy equivalence $C_*(\ltm) \cong CH^*(C^*(M), C^*(M))$.

\med
To  begin, notice that by the definitions \ref{fk},   the following diagrams commute:
$$
\begin{CD}
\Delta^k \times LM @>f_k >>  M^{k+1}  \\
@V e VV   @VV p_{ 1} V \\
M  @>> = >  M
\end{CD}
$$
where the left hand vertical map is the evaluation, $e((t_1, \cdots , t_k); \gamma) = \gamma (0)$,
and the right hand vertical map is the projection onto the first coordinate.   Pulling back the virtual
bundle $-T(M)$ defines a map of virtual bundles
$$
(f_k)_* : e^*(-TM) \la p_{1}^*(-TM), 
$$
and therefore maps of Thom spectra,  (which by abuse of notation we still call $f_k$)
 
\begin{equation}\label{tfk}
f_k : (\Delta_k)_+ \wedge \ltm \la  \mtm  \wedge (M^k)_+ .
\end{equation}
By taking adjoints, we get a map of
spectra,
$$
\begin{CD}
f : \ltm  @>\prod_k f_k >> \prod_k Map ((\Delta_k)_+ ; \ \mtm  \wedge (M^k)_+)
\end{CD}
$$
where on the right hand side the mapping spaces are maps of unital spectra.  This map is just the
induced map of Thom spectra of the map $f : LX \la \prod_{k \geq 0} Map (\Delta^k, X^{k+1})$
described in theorem \ref{jones}.  The following result is induced by theorem \ref{jones} by
passing to Thom spectra. 

Let  $\tbx_*$ be the cosimplicial spectrum defined to be the cosimplicial Thom spectrum
of  the cosimplicial virtual bundle $-TM$.  That is, the virtual bundle over the $k$ simplices  
$  M^{k+1}$  is $p_{ 1}^*(-TM)$.  Said more explicitly, $\tbx_*$ is the cosimplicial
spectrum whose $k$ - simplices are the spectrum
$$
\tbx_k = \mtm  \wedge (M^k)_+.
$$
To describe the  coface and codegeneracy maps, consider the maps
$$
\mu_L : \mtm \to M_+ \wedge \mtm \quad \text{and} \quad \mu_R : \mtm \to \mtm \wedge M_+
$$
of Thom spectra induced by the diagonal map $\Delta : M \to M \times M$. $\mu_L$  and $\mu_R$
are  the maps of Thom spectra induced by the maps of virtual bundles $\Delta_* : -TM \to
p_L^*(-TM)$ and $\Delta_* : -TM \to p_R^*(-TM)$, where $p_L$ and $p_R$ are the projection
maps $M\times M \to M$  onto the left and right coordinates respectively. We then have the
following formulas for the coface and codegeneracy maps:
 
\begin{align}\label{tcoface}
\delta_0(u; x_1, \cdots , x_{k-1}) &= (v_R; y_R, x_1, \cdots , x_{k-1}) \\
\delta_i(u; x_1 , \cdots , x_{k-1} )  &= (u; x_1, \cdots, x_{i-1}, x_i, x_i, x_{i+1}, \cdots , x_{k-1} ),
\quad  1 \leq i \leq  k-1 \notag \\
\delta_{k}(u; x_1, \cdots, x_{k-1} ) &= (v_L; x_1, \cdots , x_{k-1},  y_L ), \notag  
 \end{align}
where   $\mu_R(u) = (v_R, y_R)$,  $\mu_L(u) = (y_L, v_L)$ and
$$
\sigma_i(u; x_1, \cdots , x_{k+1} ) = (u;  x_1, \cdots, x_i, x_{i+2}, \cdots, x_{k+1} ), \quad 0 \leq i
\leq k
$$

\med
The following result is simply the application of the Thom spectrum functor for the virtual bundle
$-TM$  to theorem
\ref{jones}.

\med
\begin{theorem}\label{tjones}
Let $M$ be any closed, $d$ - dimensional manifold, and let $$\begin{CD} 
f : \ltm  @>\prod_k f_k >> \prod_k Map ((\Delta_k)_+ ; \ \mtm  \wedge (M^k)_+)
\end{CD} $$be the
product of the maps of spectra
$  f_k$, as  defined above (\ref{tfk}). Then $f$ is a homeomorphism onto its image. 
Furthermore, the image consists of sequences of maps $\{\phi_k\}$ which commute with the
coface and codegeneracy operators. 
\end{theorem}

\med
We  denote the  space of sequences of maps referred to in this theorem
$Map_{\Delta^*}(\Delta^*, \mtm  \wedge (M^k)_+)$.  This is the   total space of the
cosimplicial spectrum
$Tot(\tbx_*)$.

Notice that   the $S$ - duality between $M_+$ and $\mtm$
  defines a chain homotopy equivalence between the cochains 
$C^*(\mtm)$ (defined to be the appropriate desuspension of the cochains of the Thom space of
the normal bundle of  a fixed embedding $M \hk \br^N$ ) and the chains of the manifold
$$
C^*(\mtm) \cong C_{-*}(M_+).
$$
The maps $f_k : (\Delta_k)_+ \wedge \ltm \la  \mtm  \wedge (M^k)_+ $ then define maps of
cochains,
$$
f_k^*:  C_{-*}(M ) \otimes C^*(M) ^{\otimes k}  \cong C^*(\mtm \wedge  (M^k)_+)  \la
C^{*-k}(\ltm).
$$
Taking the dual we get a map of chain complexes

\begin{align}
(f_k)_*: C_{*-k}(\ltm)  &\la  Hom( C^*(M) ^{\otimes k} \otimes  C_{-*}(M ); \bz) \notag \\
&\cong Hom (C^*(M)^{\otimes k}; C^*(M)) \notag  \\
&= CH^k(C^*(M); C^*(M)) \notag
\end{align}

The following is then a consequence of corollary \ref{chain}, by passing to Thom spectra.

\med
\begin{corollary} \label{tchain}For any  closed manifold $M$, the chain maps $(f_k)_*$
fit together to define a 
  chain homotopy equivalence from the   chains of  Thom spectrum $\ltm$ to the Hochschild cochain
complex
$$
f_* : C_*(\ltm) \to CH^*(C^*(M); C^*(M))
$$
and so an isomorphism in homology,
$$
\begin{CD}
f_*: H_*(\ltm) @>\cong >> HH^*(C^*(M); C^*(M)).
\end{CD}
$$
\end{corollary}

\med
As mentioned in the introduction, the Hochschild cochain complex $CH^*(C^*(M); C^*(M))$ has a
cup product structure.  Namely, for any algebra $A$, if 
$$
\phi \in CH^k(A; A) = Hom (A^{\otimes k}; A) \quad \text{and} \quad  \psi \in CH^r(A; A) = Hom
(A^{\otimes r}; A),
$$
then 
$$
\phi \cup \psi \in CH^{k+r}(A; A) = Hom (A^{\otimes k+r}; A) 
$$
is defined by
$$
\phi \cup \psi (a_1 \otimes \cdots \otimes a_k \otimes a_{k+1} \otimes \cdots \otimes a_{k+r}) =
\phi(a_1 \otimes \cdots \otimes a_k)\psi(a_{k+1} \otimes \cdots \otimes a_{k+r}).
$$
For $A = C^*(M)$ (where the algebra stucture is the cup product in $ C^*(M)$),  by taking
adjoints, we can think of this as a pairing

\begin{align}\label{cup}
\cup : \left(C_*(M)^{\otimes k} \otimes C^*(M)\right) \otimes \left(C_*(M)^{\otimes r} \otimes
C^*(M)\right)  &\la C_*(M)^{\otimes k+r} \otimes C^*(M) \\
(\alpha_1 \otimes \cdots \otimes \alpha_k \otimes \theta) \otimes (\beta_1 \otimes \cdots
\otimes \beta_r \otimes \rho) &\la \alpha_1 \otimes \cdots \otimes \alpha_k \otimes \beta_1
\otimes \cdots \otimes \beta_r \otimes \theta \cup \rho. \notag  
\end{align}

Now recall that by $S$ - duality, there is a ring spectrum structure
$$
\Delta^* : \mtm \wedge \mtm \la \mtm
$$
dual to the diagonal map $\Delta : M \to M \times M$.  Passing to chains, $\Delta^* : C_*(\mtm)
\otimes C_*(\mtm) \to C_*(\mtm)$, is, with respect to the duality identification $C_*(\mtm)
\cong C^*(M)$,   the cup product on cochains
$$
\Delta^* = \cup : C^*(M) \otimes C^*(M) \to C^*(M).
$$
Thus formula  (\ref{cup}) for the cup product in Hochschild cochains is therefore realized by the
map
\begin{align}\label{mukr}
\tmu_{k,r}: \left[\mtm \wedge (M^k)_+\right]\wedge \left[\mtm \wedge (M^r)_+\right] &\la
\mtm \wedge (M^{k+r})_+ \\
(u; x_1, \cdots, x_k) \wedge (v; y_1, \cdots , y_r) &\la (\Delta^*(u,v); x_1, \cdots , x_k, y_1,
\cdots , y_r)
\notag
\end{align}

The maps $\tmu_{k,r}$ define maps of the simplices
$$
\tmu_{k,r} :  \tbx_k \wedge \tbx_r \to \tbx_{k+r}
$$
and it is straight forward to check that these maps preserve the coface and codegeneracy
operators, and induce a map of  total spectra,
$$
\tmu : Tot(\tbx_*) \wedge Tot(\tbx_*) \la Tot(\tbx_*).
$$
Notice in particular that this pairing is  $A_\infty$ associative.  This proves the following. 

\med
\begin{theorem}\label{cupcompat} Using the  homeomorphism $ f: \ltm \cong Tot(\tbx_*)$
of theorem \ref{tjones}, $\ltm$ inherits the structure of an $A_\infty$  ring spectrum,
$$
\tmu : \ltm \wedge \ltm \la \ltm
$$
which is compatible with the cup product in Hochschild cohomology.  That is, with respect to the
chain homotopy equivalence $$f_*: C_*(\ltm) \to CH^*(C^*(M); C^*(M))$$
 of corollary \ref{tchain},  the following diagram of chain complexes  commutes:
$$
\begin{CD}
C_*(\ltm) \otimes C_*(\ltm) @>\tmu >> C_*(\ltm) \\
@Vf_* \otimes f_* V \cong V    @V\cong V f_* V \\
CH^*(C^*(M); C^*(M)) \otimes CH^*(C^*(M); C^*(M))  @>>\cup > CH^*(C^*(M); C^*(M)).
\end{CD}
$$
\end{theorem}

\med
In view of theorems \ref{tjones}, \ref{tchain}, and \ref{cupcompat}, theorem \ref{cosimp} will
therefore be proven once we prove the following.

\med
\begin{theorem} Let $\mu : \ltm \wedge \ltm \to \ltm$ be the ring spectrum structure
defined in section 1.   Then the structures
$\mu $ and $\tmu : Tot(\tbx_*) \wedge  Tot(\tbx_*)  \to Tot(\tbx_*) $ are compatible in the
sense that the following diagram homotopy commutes:
$$
\begin{CD}
\ltm \wedge \ltm @>\mu >> \ltm \\
@Vf\wedge fV \simeq V   @V\simeq Vf V \\
Tot(\tbx_*) \wedge Tot(\tbx_*) @>>\tmu >  Tot(\tbx_*).  
\end{CD}
$$
\end{theorem}

\med
 \begin{proof}  Let $w_{k,r} : \Delta^{k+r} \to \Delta^k \times \Delta^r $ be the Alexander -
Whitney diagonal map.  That is, for $(x_1, \cdots , x_{k+r}) \in \Delta^{k+r}$,  then $w_{k,r}(x_1,
\cdots , x_{k+r}) = (x_1, \cdots , x_k) \times (x_{k+1}, \cdots , x_{k+r}) \in \Delta^k \times
\Delta^r$.  By the definition of the cosimplicial structure of $\tbx_*$, to prove the theorem it
suffices to prove that the following diagrams of spectra commute:

 \bfl (a)
$$
\begin{CD} 
\Delta^{k+r}_+ \wedge \ltm \wedge \ltm  @>1 \wedge \mu >>  \Delta^{k+r}_+ \wedge \ltm \\
@Vw_{k,r} \wedge 1 VV    \\
(\Delta^k \times \Delta^r)_+ \wedge \ltm \wedge \ltm  & &  @VVf_{k+r} V \\
@V f_k \wedge f_r VV   \\
\mtm \wedge (M^k)_+ \wedge \mtm \wedge (M^r)_+  @>>\tmu > \mtm \wedge (M^{k+r})_+
\end{CD}
$$

\efl 
We verify the commutativity of these diagrams in several steps.
First observe that the maps 
$$f_k \times f_r : \Delta^k \times \Delta^r \times LM \times LM \to M^{k+1} \times M^{r+1}
$$
restrict to    $\Delta^k \times \Delta^r \times LM \times_M  LM$ to define a map $f_{k,r}$ whose image
is in $M \times M^k \times M^r$ making the following diagram commute:
$$
\begin{CD}
\Delta^k \times \Delta^r \times LM \times_M LM  @>\hk >> \Delta^k \times \Delta^r \times LM
\times LM  \\
@Vf_{k,r} VV   @VVf_k \times f_r V \\
M \times M^k \times M^r  @>> \Delta > M^{k+1} \times M^{r+1}
\end{CD}
$$
where the bottom horizontal map is the diagonal map:
$$
\Delta \left(m \times (x_1, \cdots , x_k) \times (y_1, \cdots , y_r)\right) = (m, x_1, \cdots , x_k)
\times (m, y_1, \cdots, y_r).
$$
By the naturality of the  Pontrjagin - Thom construction, we therefore have a commutative diagram of spectra

\bfl
(b)
$$
\begin{CD}
(\Delta^k \times \Delta^r)_+ \wedge \ltm \wedge  \ltm @> 1 \wedge \tau >> (\Delta^k \times \Delta^r)_+
\wedge (LM \times_M LM)^{-TM} \\
@Vf_k \wedge f_r VV    @VVf_{k,r} V \\
\mtm \wedge (M^k)_+ \wedge\mtm \wedge  (M^r)_+  @>> \tau >  \mtm \wedge (M^k)_+\wedge (M^r)_+.
\end{CD}
$$

\efl
Notice further that by the definition of the maps $f_n$, $f_{k,r}$ and the loop composition
$\gamma : LM\times_M LM \to LM$ defined in the last section,  the following
diagram commutes:

$$
\begin{CD}
\Delta^{k+r} \times (LM \times_M LM)  @> 1 \times \gamma >>  \Delta^{k+r} \times LM \\
@Vw_{k,r} \times 1 VV    \\
\Delta^k \times \Delta^r \times (LM \times_M LM)   & &  @VV f_{k+r} V\\
@Vf_{k,r} VV   \\
M \times (M^k) \times (M^r)  @>> = > M^{k+r+1}.
\end{CD}
$$

Passing to Thom spectra this yields the following commutative diagram:

\bfl
(c)
$$
\begin{CD}
(\Delta^{k+r})_+ \wedge (LM \times_M LM)^{-TM}  @> 1 \wedge \gamma >>  (\Delta^{k+r})_+  \wedge \ltm \\
@Vw_{k,r} \wedge 1 VV    \\
(\Delta^k \times \Delta^r)_+ \wedge (LM \times_M LM)^{-TM}   & &  @VV f_{k+r} V\\
@Vf_{k,r} VV   \\
\mtm \wedge (M^k)_+  \wedge  (M^r)_+   @>> = > \mtm \wedge M^{k+r}_+.
\end{CD}
$$

\efl

Now consider the following diagram of spectra:

$$
\begin{CD}
(\Delta^{k+r})_+ \wedge \ltm \wedge \ltm  @> 1 \wedge \tau >> (\Delta^{k+r})_+ \wedge (LM \times_M
LM)^{-TM}  @> 1 \wedge \gamma >> (\Delta^{k+r})_+ \wedge \ltm \\
@Vw_{k,r} \wedge 1VV   @Vw_{k,r} \wedge 1VV  \\
(\Delta^k \times \Delta^r)_+ \wedge \ltm \wedge \ltm @>>1 \wedge \tau > (\Delta^k \times
\Delta^r)_+\wedge (LM \times_M LM)^{-TM} & & @VV f_{k+r} V \\
@Vf_k \wedge f_r VV   @VVf_{k,r} V \\
\mtm \wedge (M^k)_+ \wedge \mtm \wedge (M^r)_+ @>>     \tau > \mtm \wedge (M^{k+r})_+ @>> = >
\mtm \wedge (M^{k+r})_+.
\end{CD}
$$
Now the top left square in this diagram clearly commutes. The bottom left diagram is diagram (b) above,
and so it commutes.  The right hand rectangle is diagram (c) above, so it commutes.  Therefore the outside
of this diagram commutes.  Now the top horizontal composition is, by definition the map
$$
1 \wedge \mu : (\Delta^{k+r})_+ \wedge \ltm \wedge \ltm \la (\Delta^{k+r})_+ \wedge \ltm.
$$
The bottom horizontal map is seen to be  
$$
\tmu :\mtm \wedge (M^k)_+ \wedge \mtm \wedge (M^r)_+ \la \mtm \wedge (M^{k+r})_+
$$
by recalling that the ring multiplication $\Delta^* : \mtm \wedge \mtm \to \mtm$  is the  Pontrjagin  -
Thom map
$\tau : \mtm \wedge \mtm \to \mtm$ applied to the diagonal embedding  $\Delta : M \hk M \times M$. 

With these identifications, the outside of this diagram is then diagram (a) above.  As observed earlier, the
commutativity of diagram (a) proves this theorem, and this completes the proof of theorem 2.
\end{proof}

\end{document}